\documentclass[a4paper,12pt]{article}%
\usepackage{amsmath}
\usepackage{amsfonts}
\usepackage{amssymb}
\usepackage[dvips]{graphicx}%
\newtheorem{theorem}{Theorem}

\newtheorem{corollary}[theorem]{Corollary}

\newtheorem{definition}[theorem]{Definition}

\newtheorem{lemma}[theorem]{Lemma}

\newtheorem{proposition}[theorem]{Proposition}
\newtheorem{remark}[theorem]{Remark}

\newenvironment{proof}[1][Proof]{\noindent\textbf{#1.} }{\ \rule{0.5em}{0.5em}}
\begin{document}

\author{Andr\'{e} MAS$^{\ast}$, Besnik PUMO$^{\dag}\medskip$\\$\left(  ^{\ast}\right)  $ {\small Institut de Mod\'{e}lisation
Math\'{e}matique de Montpellier,}\\{\small cc 051, Place Eug\`{e}ne Bataillon,}\\{\small 34095, Montpellier Cedex 5, France,}\\{\small mas@math.univ-montp2.fr\medskip\ }\\$\left(  ^{\dag}\right)  $ {\small Unit\'{e} de Statistiques,}\\{\small UMR A 462 SAGAH,}\\{\small INH, Angers, France}}
\title{Functional linear regression with derivatives}
\date{}
\maketitle

\begin{abstract}
We introduce a new model of linear regression for random functional inputs
taking into account the first order derivative of the data. We propose an
estimation method which comes down to solving a special linear inverse
problem. Our procedure tackles the problem through a double and synchronized
penalization. An asymptotic expansion of the mean square prevision error is
given. The model and the method are applied to a benchmark dataset of
spectrometric curves and compared with other functional models.

\end{abstract}

\bigskip

\section{Introduction}

Functional Data Analysis is a well-known area of modern statistics. Advances
in computer sciences make it now possible to collect data from an underlying
continuous-time processe, say $\left(  \xi_{t}\right)  _{t\geq0}$, at high
frequencies. The traditional point of view consisting in discretizing $\left(
\xi_{t}\right)  $ at $t_{1},...,t_{p}$ and studying it by classical
multidimensional tools is outperformed by interpolation methods (such as
splines or wavelets). These techniques provide the statistician with a
reconstructed curve on which inference may be carried out through what we may
call "functional models" i.e. versions of the classical multidimensional
models designed and suited for data that are curves. Thus, functional PCA,
ANOVA or Canonical Analysis -even density estimation for curves or processes
have been investigated. We refer to Ramsay, Silverman (1997, 2002), Bosq
(2000), Ferraty Vieu (2006) for monographs on functional data analysis.
Recently many authors focused on various versions of the regression model
introduced by Ramsay and Dalzell (1991) :%
\begin{equation}
y_{i}=\int_{0}^{T}X_{i}\left(  t\right)  \rho\left(  t\right)  dt+\varepsilon
_{i} \label{mod1}%
\end{equation}
where we assume that the sample $\left(  \left(  y_{1,}X_{1}\right)
,...,\left(  y_{n},X_{n}\right)  \right)  $ is made of independent copies from
$\left(  y_{,}X\right)  .$ Each $X_{i}=\left(  X_{i}\left(  t\right)  \right)
_{t\in\left[  0,T\right]  }$ is a curve defined on the set $\left[
0,T\right]  ,$ $T>0$, $y_{i}$ is a real number, $\varepsilon_{i}$ is a white
noise and $\rho$ is an unknown function to be estimated. In other words the
$X_{i}$'s are random elements defined on an abstract probability space and
taking values in a function space, say $\mathcal{F}.$ The vector space
$\mathcal{F}$ endowed with norm $\left\Vert \cdot\right\Vert _{\mathcal{F}}$
will be described soon.We refer for instance to Cardot, Mas, Sarda (2006) or
Cai, Hall (2006) for recent results.

In this article we study a new (linear) regression model defined below derived
from (\ref{mod1}) and echoing the recent paper of Mas and Pumo (2006). The key
idea relies on the fact that most statisticians dealing with functional data
do not fully enjoy their functional properties. For instance in several models
integrals such as%
\[
\int X_{i}\left(  s\right)  X_{j}\left(  s\right)  ds
\]
are computed. The integral above is nothing but a scalar product. Nevertheless
derivatives were not given the same interest. Explicit calculations of
derivatives sometimes appear indirectly in kernel methods (when estimating the
derivatives of the density or the regression function) or through seminorms or
norms on $\mathcal{F}$. But surprisingly $X_{i}^{\prime}$ (or $X_{i}^{\left(
m\right)  }$) never appear in the models themselves whereas people dealing
with functional data often say that "derivatives contain much information,
sometimes more than the initial curves themselves". Our starting idea is the
following. Since in a functional data framework, the curve-data are
explicitely known and not just discretized, their derivatives may also be
explicitely computed. As a consequence these derivatives may be "injected" in
the model, which may enhance its prediction power. The reader is referred to
the forthcoming display (\ref{modele}) for an immediate illustration and to
Mas, Pumo (2006) for a first article dealing with a functional autoregressive
model including derivatives.

The paper is rather theoretic even if it is illustrated by a real case study.
It is organized as follows. The next section provides the mathematical
material, dealing with Hilbert spaces and linear operators, then the model is
introduced. The next section is devoted to presenting the estimation method
and its stumbling stones. The main results are given before we focus on a real
case application to food industry. The last section contains the derivation of
the theorems.

\section{About Hilbert spaces and linear operators}

Silverman (1996) provided a theoretical framework for a smoothed PCA. Jim
Ramsay (2000) enlightened the very wide scope of differential equations in
statistical modelling. Our work is in a way based on this mathematically
involved article. We are aiming at proving that derivatives may be handled in
statistical models quite easily when the space $\mathcal{F}$ is well-chosen.

The choice of the space $\mathcal{F}$ is crucial. We have to think that if
$X\in\mathcal{F}$, $X^{\prime}$ does not necessarily belong to $\mathcal{F}$
but to another space $\mathcal{F}^{\prime}$ that may be tremendously different
(larger) than $\mathcal{F}$. We decide to take $\mathcal{F}=W^{2,1},$ the
Sobolev space of order $\left(  2,1\right)  $ defined by
\[
W^{2,1}=\left\{  u\in L^{2}\left[  0,1\right]  ,u^{\prime}\in L^{2}\left[
0,1\right]  \right\}
\]
for at least three reasons :

\begin{itemize}
\item If $X\in\mathcal{F}$, $X^{\prime}\in L^{2}\left[  0,1\right]  $ which is
a well known space.

\item Both spaces are Hilbert spaces as well as%
\[
W^{2,p}=\left\{  u\in L^{2}\left[  0,1\right]  ,u^{\left(  p\right)  }\in
L^{2}\left[  0,1\right]  \right\}  .
\]
This is of great interest for mathematical reasons : bases are denumerable,
projections operators are easy to handle, covariance operators admit spectral
representations, etc.

\item The classical interpolation methods mentioned above (splines and
wavelets) provide estimates belonging to Sobolev spaces. So from a practical
point of view $W^{2,1}$ -and in general $W^{m,p},$ $\left(  m,p\right)
\in\mathbb{N}^{2},$ (see Adams and Fournier (2003) for definitions)- is a
natural space in which our curves should be imbedded.
\end{itemize}

In the sequel $W^{2,1}$ will be denoted $W$ and $W^{2,0}=L^{2}$ will be
denoted $L$ for the sake of simplicity. We keep in mind that $W$ (resp. $L$)
could be replaced by a space of higher smoothness index : $W^{2,p}$ where
$p>1$ (resp. $W^{2,p-1}$). The spaces $W$ and $L$ are separable Hilbert spaces
endowed with scalar product :%
\begin{align*}
\left\langle u,v\right\rangle _{W}  &  =\int_{0}^{1}u\left(  t\right)
v\left(  t\right)  dt+\int_{0}^{1}u^{\prime}\left(  t\right)  v^{\prime
}\left(  t\right)  dt.\\
\left\langle u,v\right\rangle _{L}  &  =\int_{0}^{1}u\left(  t\right)
v\left(  t\right)  dt
\end{align*}
and with associated norms $\left\Vert \cdot\right\Vert _{W}$ and $\left\Vert
\cdot\right\Vert _{L}$. We refer to Ziemer (1989) or to Adams and Fournier
(2003) for monographs dedicated to Sobolev spaces. Obviously if we set
$Du=u^{\prime}$ then $D$ maps $W$ onto $L$ ($D$ is the ordinary differential
operator)$.$ Furthermore Sobolev's imbedding theorem ensures that (see Adams
and Fournier (2003) Theorem 4.12 p.85) that%
\[
\left\Vert Du\right\Vert _{L}\leq C\left\Vert u\right\Vert _{W}%
\]
(where $C$ is some constant which does not depend on $u$) i.e. $D$\textbf{ }is
a bounded operator from\textbf{ }$W$\textbf{ }to\textbf{ }$L.$ This is a
crucial point to keep in mind and the fourth reason why the functional space
was chosen to be $W^{2,1}$ : the differential operator $D$ may be viewed as a
continuous linear mapping from $W$ to $L$.

Within all the paper and especially all along the proofs we will need basic
notions about operator theory. We recall a few important facts. A linear
mapping $T$ from a Hilbert space $H$ to another Hilbert space $H^{\prime}$ is
continuous whenever%
\begin{equation}
\left\Vert T\right\Vert _{\infty}=\sup_{x\in H}\frac{\left\Vert Tx\right\Vert
_{H^{\prime}}}{\left\Vert x\right\Vert _{H}}<+\infty. \label{normsup}%
\end{equation}
The adjoint of operator $T$ will be classically denoted $T^{\ast}$. Some
finite rank operators are defined by means of the tensor product : if $u$ and
$v$ belong to $H$ and $H^{\prime}$ respectively $u\otimes_{H}v$ is the
operator defined on $H$ by, for all $h\in H$~:%
\[
\left(  u\otimes_{H}v\right)  \left(  h\right)  =\left\langle u,h\right\rangle
_{H}v.
\]

\textbf{Compact operators} : Amongst linear operators the class of compact
operators is one of the best known. Compact operators generalize matrix to the
infinite-dimensional setting and feature nice properties. The general
definition of compact operators may be found in Dunford Schwartz (1988) or
Gohberg, Goldberg and Kaashoek (1991) for instance. By $\mathcal{C}_{H}$
(resp. $\mathcal{C}_{HH^{\prime}}$) we denote the space of compact operators
on the Hilbert space $H$ (resp. mapping the Hilbert space $H$ onto $H^{\prime
}$). If $T$ is a compact operator from a Hilbert space $H_{1}$ to another
Hilbert space $H_{2},$ $T$ admits the Schmidt decomposition :%
\begin{equation}
T=\sum_{k\in\mathbb{N}}\mu_{k}\left(  u_{k}\otimes v_{k}\right)
\label{Decomp.Schmidt}%
\end{equation}
where $u_{k}$ (resp. $v_{k}$) is a complete orthonormal system in $H_{1}$
(resp. in $H_{2}$) and $\mu_{k}$ are the characteristic numbers of $T$ (i.e.
the square root of the eigenvalues of $T^{\ast}T$) and
\[
\lim_{k\rightarrow+\infty}\mu_{k}=0.
\]
From (\ref{normsup}) we obtain%
\[
\left\Vert T\right\Vert _{\infty}=\sup_{k}\left\{  \mu_{k}\right\}  .
\]
When $T$ is symmetric $\mu_{k}$ is the $k^{th}$ eigenvalue of $T$ (then
$u_{k}=v_{k}$). In this situation and from (\ref{Decomp.Schmidt}) one may
define the square root of $T$ whenever $T$ maps $H$ ont $H$ and is positive :
$T^{1/2}$ is still a linear operator defined by :%
\begin{equation}
T^{1/2}=\sum_{k\in\mathbb{N}}\sqrt{\mu_{k}}\left(  u_{k}\otimes u_{k}\right)
. \label{racop}%
\end{equation}
Note that finite rank operators are always compact.\newline%
\textbf{Hilbert-Schmidt operators} : We also mention the celebrated space of
Hilbert-Schmidt operators $\mathcal{HS}\left(  H_{1},H_{2}\right)  $ - a
subspace of $\mathcal{C}_{H_{1}H_{2}}.$ Let $\left(  u_{i}\right)  _{i\geq0}$
be a basis of $H_{1}$ then $T\in\mathcal{HS}\left(  H_{1},H_{2}\right)  $
whenever%
\[
\sum_{i=1}^{+\infty}\left\Vert T\left(  u_{i}\right)  \right\Vert _{H_{2}}%
^{2}<+\infty.
\]
The space $\mathcal{HS}$ is itself a separable Hilbert space endowed with
scalar product%
\[
\left\langle T,S\right\rangle _{\mathcal{HS}}=\sum_{i=1}^{+\infty}\left\langle
T\left(  u_{i}\right)  ,S\left(  u_{i}\right)  \right\rangle _{H_{2}}%
\]
and $\left\langle T,S\right\rangle _{\mathcal{HS}}$ does not depend on the
choice of the basis $\left(  u_{i}\right)  _{i\geq0}.$ Finally the following
bound is valid for all $T\in\mathcal{HS}$ :%
\[
\left\Vert T\right\Vert _{\infty}\leq\left\Vert T\right\Vert _{\mathcal{HS}}.
\]
\newline\textbf{Unbounded operators} : If $T$ is a one to one (injective)
selfadjoint compact operator mapping a Hilbert space $H$ onto $H$, $T$ admits
an inverse $T^{-1}$. The operator $T^{-1}$ is defined on a dense (and
distinct) subspace of $H$ :%
\[
\mathcal{D}\left(  T^{-1}\right)  =\left\{  x=\sum_{p\in\mathbb{N}}x_{p}%
u_{p}:\sum_{p\in\mathbb{N}}\dfrac{x_{p}^{2}}{\mu_{p}^{2}}<+\infty\right\}  .
\]
It is unbounded which also means that $T^{-1}$ is continuous at no point for
which it is defined and $\left\Vert T^{-1}\right\Vert _{\infty}=+\infty$.

\section{The model}

We are now in position to introduce this (random input - linear) regression
model :%

\begin{equation}
y_{i}=\left\langle \phi,X_{i}\right\rangle _{W}+\left\langle \psi
,X_{i}^{\prime}\right\rangle _{L}+\varepsilon_{i} \label{modele}%
\end{equation}
where all random variables are assumed to be centered. The main result of the
paper (see next section) gives an asymptotic expansion for the mean square
prediction error in (\ref{modele}).

The unknown functions $\phi$ and $\psi$ belong to $W$ and $L$ respectively.

Obviously we are going to face two issues :

\begin{itemize}
\item Studying the identifiability of $\phi$ and $\psi$ in the model above.

\item Providing a consistent estimation procedure for $\phi$ and $\psi$.
\end{itemize}

From now on we suppose that :%
\[
\mathbf{A1}:\left\Vert X\right\Vert _{W}<M\quad a.s.
\]

This assumption could be relaxed for milder moment assumptions. We claim that
our main result holds whenever%
\[
\mathbf{A}^{\prime}\mathbf{1}:\mathbb{E}\left\Vert X\right\Vert _{W}^{8}<M.
\]
is true. But considering $\mathbf{A}^{\prime}\mathbf{1}$ would lead us to
longer and more intricate methods of proof.

\section{Estimation procedure\label{section.estimation.procedure}}

\subsection{The moment method}

Inference is based on moment formulas. From (\ref{modele}) we derive the two
following normal equation -multiply with $\left\langle X_{i},\cdot
\right\rangle $ and $\left\langle X_{i}^{\prime},\cdot\right\rangle $
successively then take expectation :%
\begin{equation}
\left\{
\begin{tabular}
[c]{l}%
$\delta=\Gamma\phi+\Gamma^{\prime}\psi,$\\
$\delta^{\prime}=\Gamma^{\prime\ast}\phi+\Gamma^{\prime\prime}\psi.$%
\end{tabular}
\ \ \ \ \ \ \ \ \ \ \right.  \label{syst}%
\end{equation}
where $\Gamma,$ $\Gamma^{\prime},$ $\Gamma^{\prime\ast},$ $\Gamma
^{\prime\prime}$ are the covariance and cross-covariance of the couple
$\left(  X_{i},X_{i}^{\prime}\right)  _{1\leq i\leq n}$ defined by :%
\begin{align*}
\Gamma &  =\mathbb{E}\left(  X\otimes_{W}X\right)  ,\;\Gamma^{\prime\ast
}=\mathbb{E}\left(  X\otimes_{W}X^{\prime}\right)  ,\\
\Gamma^{\prime}  &  =\mathbb{E}\left(  X^{\prime}\otimes_{L}X\right)
,\;\Gamma^{\prime\prime}=\mathbb{E}\left(  X^{\prime}\otimes_{L}X^{\prime
}\right)  ,
\end{align*}
and%
\[
\delta=\mathbb{E}\left(  yX\right)  \in W,\quad\delta^{\prime}=\mathbb{E}%
\left(  yX^{\prime}\right)  \in L.
\]

Under assumption $\mathbf{A1}$ or $\mathbf{A}^{\prime}\mathbf{1}$ the
covariance operators belong to $\mathcal{HS}\left(  W\right)  $,
$\mathcal{HS}\left(  W,L\right)  $, $\mathcal{HS}\left(  L,W\right)  $ or to
$\mathcal{HS}\left(  L\right)  $. Besides the covariance and cross-covariance
mentioned above are linked through the relation%
\[
\Gamma^{\prime\ast}=D\Gamma,\quad\Gamma^{\prime\prime}=D\Gamma^{\prime}.
\]

Resolving the system (\ref{syst}) is apparently easy but we should be aware of
two facts :

\begin{itemize}
\item Operators (here, $\Gamma,\Gamma^{\prime}..$.) do not commute !

\item The inverse operators of $\Gamma,$ and $\Gamma^{\prime\prime}$ do not
necessarily exist and when they do, they are unbounded, i.e. not continuous
(recall that $\Gamma,$ and $\Gamma^{\prime\prime}$ are compact operators and
that compact operators have no bounded inverses).
\end{itemize}

Before trying to solve (\ref{syst}) we will first study identifiability of the
unknown infinite dimensional parameter $\left(  \phi,\psi\right)  \in W\times
L$ in the next subsection. We complete our definitions and notations first.

We start from a sample $\left(  X_{i},X_{i}^{\prime}\right)  _{1\leq i\leq n}%
$. By $\Gamma_{n},\Gamma_{n}^{\prime},\Gamma_{n}^{\prime\ast},\Gamma
_{n}^{\prime\prime},\delta_{n}$ and $\delta_{n}^{\prime}$ we denote the
empirical counterparts of the operators and vectors introduced above and based
on the sample $\left(  y_{i},X_{i},X_{i}^{\prime}\right)  _{1\leq i\leq n}$.
For example :%
\begin{align}
\Gamma_{n}  &  =\dfrac{1}{n}\sum_{k=1}^{n}X_{k}\otimes_{W}X_{k},\label{esti}\\
\Gamma_{n}^{\prime}  &  =\dfrac{1}{n}\sum_{k=1}^{n}X_{k}^{\prime}\otimes
_{L}X_{k},\nonumber\\
\delta_{n}  &  =\dfrac{1}{n-1}\sum_{k=1}^{n-1}y_{k}X_{k}.\nonumber
\end{align}

\subsection{Identifiability}

Both equations in (\ref{syst}) are the starting point of the estimation
procedure. We should make sure that solutions to these equations are well and
uniquely defined. Suppose for instance that \textrm{Ker}$\Gamma\neq\left\{
0\right\}  $ and take $h$ in it. Now set $\widetilde{\phi}=\phi+h.$ Then
\[
\Gamma\widetilde{\phi}=\Gamma\phi+\Gamma h=\Gamma\phi.
\]
So $\Gamma\widetilde{\phi}=\Gamma\phi$ and since $\Gamma^{\prime\ast}=D\Gamma$
it is plain that $\Gamma^{\prime\ast}\widetilde{\phi}=\Gamma^{\prime\ast}%
\phi.$ Consequently $\widetilde{\phi}$ is another solution to (\ref{syst}).
There are indeed even infinitely many solutions in the space $\phi
+$\textrm{Ker}$\Gamma$. For similar reasons about $\psi$ we should impose
$\mathrm{Ker}T=\left\{  0\right\}  $ for $T=\left\{  \Gamma,\Gamma^{\prime
},\Gamma^{\prime\ast},\Gamma^{\prime\prime}\right\}  .$ It turns out that the
only necessary assumption is%
\[
\mathbf{A2}:\mathrm{Ker}\Gamma=\mathrm{Ker}\Gamma^{\prime\prime}=\left\{
0\right\}  \mathrm{.}%
\]
It is easily seen that $\mathbf{A2}$ implies $\mathrm{Ker}\Gamma^{\prime
}=\mathrm{Ker}\Gamma^{\prime\ast}=\left\{  0\right\}  .$ With other words we
suppose that both operators $\Gamma$ and $\Gamma^{\prime\prime}$ above are one
to one.

We are now ready to solve the identification problem.

\begin{proposition}
\label{Ident}The couple $\left(  \phi,\psi\right)  \in W\times L$ is
identifiable for the moment method proposed in (\ref{syst}) if and only if
$\mathbf{A2}$ holds and $\left(  \phi,\psi\right)  \notin\mathcal{N}$ where
$\mathcal{N}$ is the vector subspace of $W\times L$ defined by :%
\begin{equation}
\mathcal{N=}\left\{  \left(  \phi,\psi\right)  :\phi+D^{\ast}\psi=0\right\}  .
\label{fluke}%
\end{equation}

\end{proposition}

The above Proposition is slightly abstract but (\ref{fluke}) may be simply
rewritten: $\left(  \phi,\psi\right)  \in\mathcal{N}$ whenever for all
function $f$ in $W,$%
\[
\int\left(  f\phi+f^{\prime}\phi^{\prime}+f^{\prime}\phi\right)  =0
\]

Note that $\mathcal{N}$ is a closed set in $W\times L$. From now on we will
assume that :%
\[
\mathbf{A3}:\left(  \phi,\psi\right)  \notin\mathcal{N}\text{.}%
\]

\section{Definition of the estimates}

The estimates stem from (\ref{syst}) which is a non invertible system. Under
assumption $\mathbf{A2}$ the solution exists and is unique :%
\begin{equation}
\left\{
\begin{tabular}
[c]{l}%
$\phi=\left(  \Gamma-\Gamma^{\prime}\Gamma^{\prime\prime-1}\Gamma^{\prime\ast
}\right)  ^{-1}\left[  \delta-\Gamma^{\prime}\Gamma^{\prime\prime-1}%
\delta^{\prime}\right]  ,$\\
$\psi=\left(  \Gamma^{\prime\prime}-\Gamma^{\prime\ast}\Gamma^{-1}%
\Gamma^{\prime}\right)  ^{-1}\left[  \delta^{\prime}-\Gamma^{\prime\ast}%
\Gamma^{-1}\delta\right]  .$%
\end{tabular}
\ \ \ \ \ \ \ \ \right.  \label{soluce}%
\end{equation}
Let us denote
\begin{align*}
S_{\phi}  &  =\Gamma-\Gamma^{\prime}\Gamma^{\prime\prime-1}\Gamma^{\prime\ast
},\\
S_{\psi}  &  =\Gamma^{\prime\prime}-\Gamma^{\prime\ast}\Gamma^{-1}%
\Gamma^{\prime}.
\end{align*}
The reader should note two crucial facts. On the one hand $\Gamma^{-1}$ and
$\Gamma^{\prime\prime-1}$ are unbouded operators but closed graphs argument
ensure that $\Gamma^{\prime}\Gamma^{\prime\prime-1}\delta^{\prime}$ and
$\Gamma^{\prime\ast}\Gamma^{-1}\delta$ exist in $W$ and $L$ respectively. On
the other hand $\delta-\Gamma^{\prime}\Gamma^{\prime\prime-1}\delta^{\prime}$
(resp. $\delta^{\prime}-\Gamma^{\prime\ast}\Gamma^{-1}\delta$) belong to the
domain of the unbounded operator $S_{\phi}^{-1}$ (resp. $S_{\psi}^{-1}$) which
also ensures the finiteness of both solutions given in the display above.

Finding approximations to the solutions of (\ref{soluce}) is known in the
mathematical literature as "solving a linear inverse problem". The book by
Tikhonov and Arsenin (1977) -as many other references therein- is devoted to
this theory well-known in image reconstruction. The unboundedness of $S_{\phi
}^{-1}$ may cause large variation of $S_{\phi}^{-1}x$ even for small
variations of $x$. This lack of stability turns out to damage, as well as the
traditional "curse of dimensionality", the rates of convergence of our estimates.

Unfortunately we cannot simply replace "theoretical" operators and vectors by
their empirical estimates because $\Gamma_{n}$ and $\Gamma_{n}^{\prime\prime}$
are not invertible. Indeed they are finite-rank operators (for example the
image of $\Gamma$ is \textrm{span}$\left(  X_{1},...,X_{n}\right)  $) hence
not even injective. We are classically going to add a small perturbation to
regularize $\Gamma_{n}$ and $\Gamma_{n}^{\prime\prime}$ (see Tikhonov and
Arsenin (1977)) and another one for $S_{\phi}^{-1}$ and make them invertible.
At last $\Gamma^{-1}$ is approximated by $\Gamma_{n}^{\dag}=\left(  \Gamma
_{n}+\alpha_{n}I\right)  ^{-1},$ $\Gamma^{\prime\prime-1}$ by $\Gamma
_{n}^{\prime\prime\dag}=\left(  \Gamma_{n}^{\prime\prime}+\alpha_{n}I\right)
^{-1}$ and $S_{\phi}^{-1}$ by $\left(  S_{n,\phi}+\beta_{n}I\right)  ^{-1}$
where%
\[
S_{n,\phi}=\Gamma_{n}-\Gamma_{n}^{\prime}\left(  \Gamma_{n}^{\prime
\prime\dagger}\right)  \Gamma_{n}^{\prime\ast}.
\]
and $\alpha_{n}>0,$ $\beta_{n}>0.$ We also set :%
\begin{align}
S_{n,\psi}  &  =\Gamma_{n}^{\prime\prime}-\Gamma_{n}^{\prime\ast}\left(
\Gamma_{n}^{\dagger}\right)  \Gamma_{n}^{\prime},\label{Snpsi}\\
u_{n,\phi}  &  =\delta_{n}-\Gamma_{n}^{\prime}\left(  \Gamma_{n}^{\prime
\prime\dagger}\right)  \delta_{n}^{\prime},\label{Tnphi}\\
u_{n,\psi}  &  =\delta_{n}^{\prime}-\Gamma_{n}^{\prime\ast}\left(  \Gamma
_{n}^{\dagger}\right)  \delta_{n}. \label{Tnpsi}%
\end{align}

In the sequel we will assume that both strictly positive sequences $\alpha
_{n}$ and $\beta_{n}$ decay to zero in order to get the asymptotic convergence
of the estimates.

\begin{definition}
The estimate of the couple $\left(  \phi,\psi\right)  $ is $\left(
\widehat{\phi}_{n},\widehat{\psi}_{n}\right)  $ based on (\ref{soluce}) and
defined by :%
\begin{equation}
\left\{
\begin{tabular}
[c]{l}%
$\widehat{\phi}_{n}=\left(  S_{n,\phi}+\beta_{n}I\right)  ^{-1}u_{n,\phi
},\smallskip$\\
$\widehat{\psi}_{n}=\left(  S_{n,\psi}+\beta_{n}I\right)  ^{-1}u_{n,\psi}.$%
\end{tabular}
\ \right.  \label{estim}%
\end{equation}
The predictor is defined as%
\[
\widehat{y}_{n+1}=\left\langle \widehat{\phi}_{n},X_{n+1}\right\rangle
_{W}+\left\langle \widehat{\psi}_{n},X_{n+1}^{\prime}\right\rangle _{L}.
\]

\end{definition}

\section{Main results and comments}

In Mas, Pumo (2006) the authors obtained convergence in probability for their
estimates in a quite different model. We are now in position to assess deeper
results. Mean square prediction error is indeed given an asymptotic
development depending on both smoothing sequences $\alpha_{n}$ and $\beta_{n}$.

Before stating the main result of this article, we give and comment the next
and last assumption :%
\begin{equation}
\mathbf{A4}:\left\{
\begin{tabular}
[c]{l}%
$\left\Vert \Gamma^{-1/2}\phi\right\Vert _{W}<+\infty$\\
$\left\Vert \left(  \Gamma^{\prime\prime}\right)  ^{-1/2}\psi\right\Vert
_{L}<+\infty$%
\end{tabular}
\ \ \ \ \ \right.  \label{A4}%
\end{equation}

For the definition of $\Gamma^{-1/2}$ and $\Gamma^{\prime\prime-1/2}$ we refer
to (\ref{racop}). Let us explain briefly what both conditions in (\ref{A4})
mean. To that aim we rewrite the first by developing $\Gamma^{-1/2}\phi$ in a
basis of eigenvectors of $\Gamma,$ say $u_{p}$%
\[
\Gamma^{-1/2}\phi=\sum_{p=1}^{+\infty}\dfrac{\left\langle \phi,u_{p}%
\right\rangle }{\sqrt{\lambda_{p}}}u_{p}%
\]
hence%
\[
\left\Vert \Gamma^{-1/2}\phi\right\Vert _{W}^{2}=\sum_{p=1}^{+\infty}%
\dfrac{\left\langle \phi,u_{p}\right\rangle ^{2}}{\lambda_{p}}%
\]
The first part of assumption $\mathbf{A4}$ tells us that "$\left\langle
\phi,u_{p}\right\rangle $ should tend to zero quickly enough with respect to
$\lambda_{p}$". In other words $\phi$ should belong to an ellipso\"{\i}d of
$W$ which may be more or less "flat" depending on the rate of decay of the
$\lambda_{p}$'s to zero. Assumption $\mathbf{A4}$ is in fact \textbf{a
regularity condition} on functions $\phi$ and $\psi$ : function $\phi$ (resp.
$\psi$) should be smoother than $X$ (resp. $X^{\prime}$).\bigskip

We could try and state convergence results for $\widehat{\phi}_{n}$ and
$\widehat{\psi}_{n}$ separatedly but it turns out that :

\begin{itemize}
\item The real statistical interest of the model relies on its predictive
power. The statistician is mainly interested in $\widehat{y}_{n+1},$ not in
$\widehat{\phi}_{n}$ and $\widehat{\psi}_{n}$ in a first attempt. The issue of
goodness of fit tests (involving $\phi$ and $\psi$ alone) is beyond the scope
of this article.

\item Considering the mean square norm of $\left\langle \widehat{\phi}%
_{n},X_{n+1}\right\rangle _{W}$ (instead of $\widehat{\phi}_{n}$ or even of
$\left\langle \widehat{\phi}_{n},x\right\rangle _{W}$ for a nonrandom $x$) has
a smoothing effect on our estimates and partially counterbalance the side
effects of the underlying inverse problem as will be seen within the proofs
(especially along Lemma \ref{placebo}).
\end{itemize}

Turning to $\widehat{y}_{n+1},$ the next question is : what should we compare
$\widehat{y}_{n+1}$ with~? The right answer is not $y_{n+1}.$ Obviously we
could, but it is also plain that, due to the random $\varepsilon_{n+1}$ the
best possible prediction for $y_{n+1}$ knowing $X_{n+1}$ (or even the "past"
i.e. $X_{1},...,X_{n}$) is the conditional expectation :%
\[
y_{n+1}^{\ast}=\mathbb{E}\left(  y_{n+1}|X_{1},...,X_{n+1}\right)
=\left\langle \phi,X_{n+1}\right\rangle _{W}+\left\langle \psi,X_{n+1}%
^{\prime}\right\rangle .
\]
We are now ready to state the main theoretical result of this article.

\begin{theorem}
\label{conv}When assumptions $\mathbf{A1-A4}$ hold the following expansion is
valid for the prediction mean square error :%
\[
\mathbb{E}\left(  \widehat{y}_{n+1}-y_{n+1}^{\ast}\right)  ^{2}=O\left(
\frac{\beta^{2}}{\alpha^{2}}\right)  +O\left(  \frac{1}{\alpha^{2}\beta^{2}%
n}\right)
\]

\end{theorem}

\begin{remark}
Replacing $y_{n+1}^{\ast}$ with $y_{n+1}$ is still possible. We may easily
prove that~:%
\[
\mathbb{E}\left(  \widehat{y}_{n+1}-y_{n+1}\right)  ^{2}=\mathbb{E}\left(
\widehat{y}_{n+1}-y_{n+1}^{\ast}\right)  ^{2}+\sigma_{\varepsilon}^{2}.
\]

\end{remark}

\begin{corollary}
From Theorem \ref{conv} above an optimal choice for $\beta$ is $\beta^{\ast
}\asymp n^{-1/4},$ then the convergence rate is :%
\[
\mathbb{E}\left(  \widehat{y}_{i}-y_{i}^{\ast}\right)  ^{2}=O\left(  \frac
{1}{\alpha^{2}n^{1/2}}\right)
\]
and may be quite close from $1/n^{1/2}$.
\end{corollary}

The proof of the Corollary will be omitted. Studying the optimality of this
rate of convergence over the classes of functions defined by $\mathbf{A4}$ is
beyond the scope of this article but could deserve more attention.

\begin{remark}
Originally the linear model (\ref{modele}) is subject to serious
multicolinearity troubles since $X_{n}^{\prime}=DX_{n}.$ Even if the curve
$X_{n}^{\prime}$ usually looks quite different from $X_{n},$ there is a total
stochastic dependence between them. The method used in this article to tackle
this problem (as well as the intrinsic "inverse problem" aspects related to
the inversion of the covariance operators $\Gamma$ and $\Gamma^{\prime\prime}%
$) is new up to the authors' knowledge. As it can be seen through above at
display (\ref{estim}) or in the proofs below, it relies on a double
penalization technique first by the index $\alpha_{n}$ then by $\beta_{n}$
linking both indexes in order to suppress the bias terms asymptotically .
\end{remark}

\section{An application to spectrometric data}

In this section we will present an application of the Functional Linear
Regression with Derivatives (FLRD) introduced in this paper to a spectroscopic
calibration problem. Quantitative NIR (near-infrared) spectroscopy is used to
analyze food and agricultural materials. The NIR spectrum of a sample is a
continuous curve giving the absorption, that is $\log_{10}1/R$ where $R$ is
the reflection of the sample, against wavelength measured in nanometers (nm).

In the \texttt{cookie} example considered here the aim is to predict the
percentage of each ingredient $y$ given the NIR spectrum $x$ of the sample
(see Osborne et al. (1984) for a full description of the experiment). The
constituents under investigation are: fat, sucrose, dry flour, and water.
There were 39 samples in the calibration set, sample number 23 having been
excluded from the original 40 as an outlier, and a further validation set with
31 samples, again after the exclusion of one outlier.

\begin{figure}
[ptb]
\begin{center}
\includegraphics[
height=3.1099in,
width=4.9035in
]%
{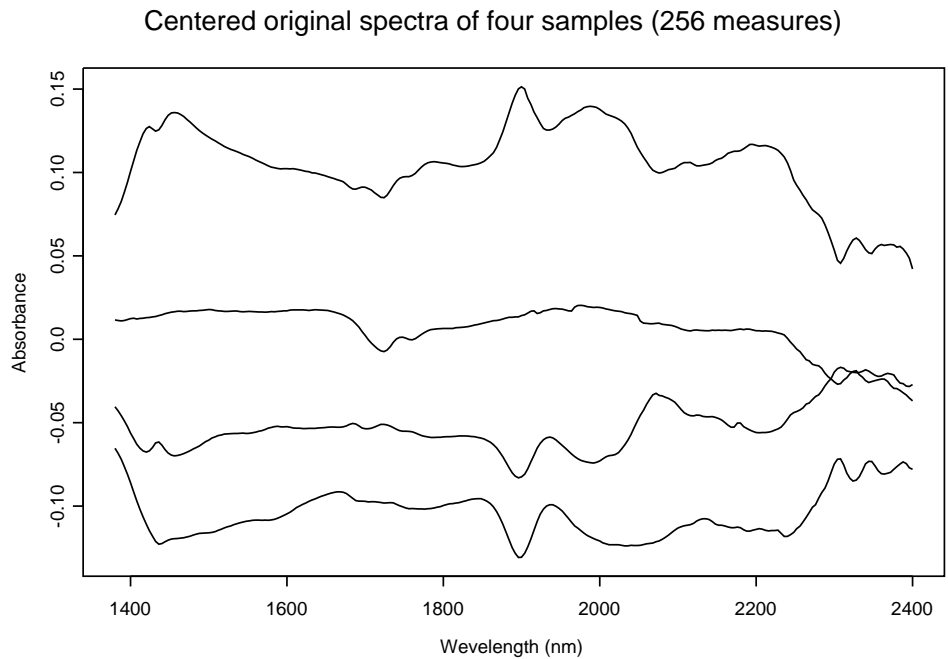}%
\caption{Centered original spectra of four samples (256 measures)}%
\label{Cookie_spectres}%
\end{center}
\end{figure}

An NIR reflectance spectrum is available for each dough. The original spectral
data consists of 700 points measured from 1100 to 1498 nm in steps of 2 nm.
Following Brown et al. (2001) we reduced the number of spectral points to 256
by considering only the spectral range 1380-2400 nm in step of 4 nm. Samples
of centered spectra are plotted in Figure \ref{Cookie_spectres}.\bigskip\


A classical tool employed in the chemiometric literature for the prediction of
$y$ knowing the associated NIR spectra $(x_{j},j=1,\ldots,256)$ is the linear
model:%
\begin{equation}
y=\sum_{j=1,256}\theta_{j}x_{j}+\epsilon\label{Mult.Lin.Mod}%
\end{equation}
The problem then is to use the calibration data to estimate the unknown
parameters $\theta_{j}$. Clearly in this application since $39\ll256$ the
ordinary least squares fails and many authors proposed to use alternative
methods to tackle the problem: principal component regression (PCR) or partial
least squares regression (PLS). We invite the reader to look at the paper of
Frank and Friedman (1993) for a statistical view of some chemiometrics
regression tools.

Following an idea of Hastie and Mallows, in their discussion of Frank and
Friedman's paper, we consider a spectrum as a functional observation. The
functional Linear Regression (FLR) corresponding to the model
\ref{Mult.Lin.Mod} defined above is:%
\[
y=\int_{\delta}x(t)\theta(t)dt+\varepsilon
\]
where $y$ is a scalar random variable, $x$ a real function defined on
$\delta=[1100,2400]$\ and $\theta(t)$ the unknown parameter function. Brown et
al. (2001), Ferraty and Vieu (2003), Marx and Eilers (2002) or Amato et al.
(2006) used such a model for a prediction problem with spectrometric data.

The model FLRD introduced in this paper can be written as:%
\[
y=\int_{\delta}x(t)\phi(t)dt+\int_{\delta}x^{\prime}(t)\psi(t)dt+\varepsilon
\]
where $\phi(t)$ and $\psi(t)$ are unknown functions (see display
(\ref{modele}) for an equivalent definition). In this paragraph we compare the
performance of PCR, PLS, FLR, FLRD, Spline Smoothing model proposed by Cardot,
Ferraty and Sarda (2006) and Bayes wavelet predictions proposed by Brown et
al. (2001).

We used the calibration data set for the estimation of parameter functions
$\phi(t)$ and $\psi(t)$ and validation data for calculation of the MSEP (Mean
Squared Error of Predictions):%
\[
MSEP=\frac{1}{31}\sum_{j=1}^{31}(y_{j}-\hat{y}_{j})^{2}%
\]
where $\hat{y}_{j}$\ is the prediction of $y_{j}$ obtained by the model with
estimated parameters. The choice of the parameters $\alpha$ and $\beta$ is
crucial for the prediction model. We used a cross-validation approach based on
the evaluation of the standard error of prediction $CVMSEP$:
\[
CVMSEP(\alpha,\beta)=\frac{1}{39}\sum_{i=1}^{39}[\frac{1}{38}\sum
\nolimits_{j=1}^{38}(y_{i}^{c}-\hat{y}_{j}^{c}(i;\alpha,\beta))^{2}],
\]
where $\hat{y}_{j}^{c}(i;\alpha,\beta)$ denotes the prediction of $y_{j}^{c}%
$\ in the calibration set without sample $i$. Results for different methods of
prediction of four ingredients are displayed in Table \ref{Table.MSEP}. We
used B-spline basis ($k=100$) for obtaining predictions with Spline Smoothing,
Spline Ridge RLF and Spline RLFD methods. For each of those methods we give
the values of the smoothing or penalty parameters based on an analogous
cross-validation approach.

\begin{table}[h]
\begin{center}%
\begin{tabular}
[c]{|l|cccc|}\hline
& \multicolumn{4}{|c|}{MSE Validation}\\
Method and parameters & Fat & Sugar & Flour & Water\\\hline
PLS & 0.151 & 0.583 & 0.375 & 0.105\\
PCR & 0.160 & 0.614 & 0.388 & 0.106\\
Spline Smoothing ($k_{n}=8$) & 0.546 & 0.471 & 2.226 & 0.183\\
Spline Ridge FLR ($\beta=0.00002$) & 0.044 & 0.494 & 0.318 & 0.087\\
Spline FLRD ($\alpha=0.07,\beta=0.15$) & 0.092 & 0.450 & 0.332 & 0.069\\
Bayes Wavelet & 0.063 & 0.449 & 0.348 & 0.050\\\hline
\end{tabular}
\end{center}
\caption{$MSEP$ criterion for all models (see Brown et al. for results of PLS,
PCR and Bayes wavelet methods).}%
\label{Table.MSEP}%
\end{table}

We note that functional approaches work better then PLS or PCR methods for the
four predicted variables with respect to $MSEP$ criterion. Our simulation, as
noted also by Marx and Eilers (2002), show that functional methods lead to
more stable prediction. The Spline FLRD method produces in general equivalent
results in terms of predictions with the best methods presented in table
\ref{Table.MSEP}.

\bigskip

\section{Proofs}

In the sequel $M$ and $M^{\prime}$ will stand for constants.

Let $S$ and $T$ be two selfadjoint linear operators on a Hilbert space $H,$ we
denote $T\ll S$ whenever for all $x$ in $H,$ $\left\langle Tx,x\right\rangle
\leq\left\langle Sx,x\right\rangle $ then $\left\Vert T\right\Vert _{\infty
}\leq\left\Vert S\right\Vert _{\infty}$.

The norm in the space $L^{2}\left(  B\right)  $ where $\left(  B,\left\Vert
\cdot\right\Vert _{B}\right)  $ is a Banach space is defined the following way
: let $X$ be a random element in the Banach space $B,$ then
\[
\left\Vert X\right\Vert _{L^{2}\left(  B\right)  }=\left(  E\left\Vert
X\right\Vert _{B}^{2}\right)  ^{1/2}%
\]
When the notation is not ambiguous we systematically drop the index $B$ i.e :
$\left\Vert X\right\Vert _{L^{2}}=\left(  E\left\Vert X\right\Vert _{B}%
^{2}\right)  ^{1/2}$.

\subsection{Preliminary facts :}

In order to gain some clarity in the proofs and to alleviate them we first
list a few results stemming from operator or probabillity theory.

\textbf{Fact 1}: If $T$ is a positive operator (either random or not),
$T+\gamma I$ is invertible for all $\gamma>0$ with bounded inverse and
$\left\Vert \left(  T+\gamma I\right)  ^{-1}\right\Vert _{\infty}\leq
\gamma^{-1}$. Hence%
\begin{equation}
\left\Vert \Gamma_{n}^{\dag}\right\Vert _{\infty}=\left\Vert \Gamma^{\dag
}\right\Vert _{\infty}=\left\Vert \Gamma_{n}^{\prime\prime\dag}\right\Vert
_{\infty}=\left\Vert \Gamma^{\prime\prime\dag}\right\Vert _{\infty}%
=\alpha^{-1} \label{norm.res}%
\end{equation}

\textbf{Fact 2}: As a consquence of assumption $\mathbf{A1}$ and of the strong
law of large numbers for Hilbert valued random elements (see Ledoux, Talagrand
(1991) Chapter 7),%
\[
T_{n}\underset{n\rightarrow+\infty}{\rightarrow}T\quad a.s.
\]
whenever $T_{n}=\Gamma_{n},\Gamma_{n}^{\prime},\Gamma_{n}^{\prime\ast}%
,\Gamma_{n}^{\prime\prime}$ (resp. $T=\Gamma,\Gamma^{\prime},\Gamma
^{\prime\ast},\Gamma^{\prime\prime}$) since all theses random operators may be
rewritten as sums of i.i.d. random variables. These sequences of random
operators are almost surely bounded%
\begin{equation}
\sup_{n}\left\Vert T_{n}\right\Vert _{\infty}\leq M\quad a.s.
\label{norm.delta}%
\end{equation}
which also means that
\begin{equation}
\max\left(  \sup_{n}\left\Vert \delta_{n}\right\Vert _{W},\sup_{n}\left\Vert
\delta_{n}^{\prime}\right\Vert _{L}\right)  \leq M^{\prime} \label{asd}%
\end{equation}
since (for instance) $\delta_{n}=\Gamma_{n}\phi+\Gamma_{n}^{\prime}\psi+e_{n}$
where $e_{n}$ is again a sum of i.i.d random elements :%
\[
e_{n}=\frac{1}{n}\sum_{k=1}^{n}X_{k}\varepsilon_{k}%
\]
We also set%
\[
e_{n}^{\prime}=\frac{1}{n}\sum_{k=1}^{n}X_{k}^{\prime}\varepsilon_{k}%
\]
(see below for details).

\textbf{Fact 3}: The Central Limit Thorem in Hilbert spaces (or standards
results on rates of convergence for Hilbert valued random elements in square
norm) provide a rate in the $L^{2}$ convergence of several random variables of
interest in the proofs. See for instance Ledoux, Talagrand (1991) or Bosq
(2000) . Whenever $T_{n}=\Gamma_{n},\Gamma_{n}^{\prime},\Gamma_{n}^{\prime
\ast},\Gamma_{n}^{\prime\prime}$ (resp. $T=\Gamma,\Gamma^{\prime}%
,\Gamma^{\prime\ast},\Gamma^{\prime\prime}$) we have $\mathbb{E}\left\Vert
T_{n}-T\right\Vert _{\mathcal{HS}}^{2}=O\left(  \frac{1}{n}\right)  $ hence
\begin{equation}
\left\Vert T_{n}-T\right\Vert _{L^{2}\left(  \mathcal{HS}\right)  }=O\left(
\dfrac{1}{\sqrt{n}}\right)  \label{L2}%
\end{equation}
since all theses random operators may be rewritten as sums of i.i.d. random variables.

We begin with proving Proposition \ref{Ident}.

\textbf{Proof of Proposition \ref{Ident} :}

The method of the proof may be adapted from the model studied in Mas, Pumo
(2006). The couple $\left(  \phi,\psi\right)  $ will be identified whenever,
for any other couple $\left(  \phi_{a},\psi_{a}\right)  $, if
\[
\left\{
\begin{tabular}
[c]{l}%
$\delta=\Gamma\phi+\Gamma^{\prime}\psi=\Gamma\phi_{a}+\Gamma^{\prime}\psi
_{a},$\\
$\delta^{\prime}=\Gamma^{\prime\ast}\phi+\Gamma^{\prime\prime}\psi
=\Gamma^{\prime\ast}\phi_{a}+\Gamma^{\prime\prime}\psi_{a}.$%
\end{tabular}
\right.
\]
$\left(  \phi_{a},\psi_{a}\right)  =\left(  \phi,\psi\right)  $. This will be
true if%
\[
\left\{
\begin{tabular}
[c]{l}%
$\Gamma\left(  \phi-\phi_{a}\right)  +\Gamma^{\prime}\left(  \psi-\psi
_{a}\right)  =0,$\\
$\Gamma^{\prime\ast}\left(  \phi-\phi_{a}\right)  +\Gamma^{\prime\prime
}\left(  \psi-\psi_{a}\right)  =0.$%
\end{tabular}
\right.
\]

This means that the couple $\left(  \phi-\phi_{a},\psi-\psi_{a}\right)  $
belongs to the kernel of the linear operator defined blockwise on $W\times L$
by :%
\[
\left(
\begin{array}
[c]{cc}%
\Gamma & \Gamma^{\prime}\\
\Gamma^{\prime\ast} & \Gamma^{\prime\prime}%
\end{array}
\right)  .
\]
As $\Gamma^{\prime\ast}=D\Gamma$ and $\Gamma^{\prime\prime}=D\Gamma^{\prime}$,
the Proposition will be proved if the blockwise operator defined on $W\times
L$ and with values in $W$ :%
\[
\left(
\begin{array}
[c]{cc}%
\Gamma & \Gamma^{\prime}%
\end{array}
\right)  =\left(
\begin{array}
[c]{cc}%
\Gamma & \Gamma D^{\ast}%
\end{array}
\right)
\]
is one to one. It is plain that the kernel of this operator is precisely the
space $\mathcal{N}$ that appears at display (\ref{fluke}).

This finishes the proof of the Proposition.

The next two general Propositions are proved for further purpose.

\begin{proposition}
\label{lodge}%
\begin{align*}
\sup_{n}\left\Vert \left(  \Gamma_{n}^{\prime\prime\dag}\right)  ^{1/2}%
\Gamma_{n}^{\prime\ast}\right\Vert _{\infty}  &  <M\quad a.s.,\\
\sup_{n}\left\Vert \Gamma_{n}^{\prime}\left(  \Gamma_{n}^{\prime\prime\dag
}\right)  ^{1/2}\right\Vert _{\infty}  &  <M\quad a.s.,\\
\sup_{n}\left\Vert \left(  \Gamma^{\prime\prime\dag}\right)  ^{1/2}%
\Gamma^{\prime\ast}\right\Vert _{\infty}  &  <M,\\
\sup_{n}\left\Vert \Gamma^{\prime}\left(  \Gamma^{\prime\prime\dag}\right)
^{1/2}\right\Vert _{\infty}  &  <M.
\end{align*}

\end{proposition}

\begin{proof}
We prove only the first bound since the method may be copied for the other
ones. Set $R_{n}=D\Gamma_{n}^{1/2}$ then :%
\begin{align*}
\Gamma_{n}^{\prime\prime}  &  =R_{n}R_{n}^{\ast},\\
\left(  \Gamma_{n}^{\prime\prime\dag}\right)  ^{1/2}\Gamma_{n}^{\prime\ast}
&  =\left(  R_{n}R_{n}^{\ast}+\alpha I\right)  ^{-1/2}R_{n}\Gamma_{n}^{1/2}.
\end{align*}
At last,%
\[
\left\Vert \left(  \Gamma_{n}^{\prime\prime\dag}\right)  ^{1/2}\Gamma
_{n}^{\prime\ast}\right\Vert _{\infty}\leq\left\Vert \left(  R_{n}R_{n}^{\ast
}+\alpha I\right)  ^{-1/2}R_{n}\right\Vert _{\infty}\left\Vert \Gamma
_{n}^{1/2}\right\Vert _{\infty}.
\]
It is plain that%
\[
\sup_{n}\left\Vert \Gamma_{n}^{1/2}\right\Vert _{\infty}\leq M\quad a.s.
\]
If the Schmidt decomposition of $R_{n}$ is :%
\[
R_{n}=\sum_{k\in\mathbb{N}}\mu_{k,n}\left(  u_{k,n}\otimes v_{k,n}\right)  ,
\]
$(u_{k,n}\in W,\ v_{k,n}\in L)$ it is simple algebra to get :%
\begin{equation}
\left(  R_{n}R_{n}^{\ast}+\alpha I\right)  ^{-1/2}R_{n}=\sum_{k\in\mathbb{N}%
}\frac{\mu_{k,n}}{\sqrt{\mu_{k,n}^{2}+\alpha}}\left(  u_{k,n}\otimes
v_{k,n}\right)  \label{mano}%
\end{equation}
which yields $\left\Vert \left(  R_{n}R_{n}^{\ast}+\alpha I\right)
^{-1/2}R_{n}\right\Vert _{\infty}=\sup_{k}\left\{  \frac{\mu_{k,n}}{\sqrt
{\mu_{k,n}^{2}+\alpha}}\right\}  \leq1.$
\end{proof}

\begin{proposition}
\label{KT}%
\[
\left\Vert \left(  S_{n}+\beta I\right)  ^{-1}\right\Vert _{\infty}\leq
\dfrac{1}{\beta}.
\]

\end{proposition}

\begin{proof}
The proof of this Lemma is similar to Lemma 7.4 in Mas, Pumo (2006). It was
then proved for $S$ instead of $S_{n}$ and all operators should be changed to
their empirical counterparts (e.g : $\Gamma_{n}$ insted of $\Gamma$). We give
a sketch of it. The proof relies on the Schmidt decomposition of $S_{n}.$ One
would get%
\[
S_{n}=\Gamma_{n}^{1/2}\Lambda_{n}\left(  \alpha\right)  \Gamma_{n}^{1/2}%
\]
where $\Lambda_{n}\left(  \alpha\right)  $ and $\Gamma_{n}^{1/2}$ are
symmetric positive operators, which implies that $S_{n}$ itself is positive.
It suffices then to apply \textbf{Fact 2 }(see the "Preliminary facts"
subsection) to get the desired result.
\end{proof}

\subsection{Outline of the proof of Theorem \ref{conv} :}

The following bound is valid :%
\begin{align*}
&  \left[  \widehat{y}_{n+1}-\left(  \left\langle \phi,X_{n+1}\right\rangle
_{W}+\left\langle \psi,X_{n+1}^{\prime}\right\rangle _{L}\right)  \right]
^{2}\\
&  =\left(  \left\langle \phi-\widehat{\phi},X_{n+1}\right\rangle
_{W}+\left\langle \psi-\widehat{\psi},X_{n+1}^{\prime}\right\rangle
_{L}\right)  ^{2}\\
&  \leq2\left[  \left\langle \phi-\widehat{\phi},X_{n+1}\right\rangle _{W}%
^{2}+\left\langle \psi-\widehat{\psi},X_{n+1}^{\prime}\right\rangle _{L}%
^{2}\right]  .
\end{align*}
Then
\begin{align*}
\mathbb{E}\left\langle \phi-\widehat{\phi},X_{n+1}\right\rangle _{W}^{2}  &
=\mathbb{E}\left[  \mathbb{E}\left\langle \phi-\widehat{\phi},X_{n+1}%
\right\rangle _{W}^{2}|X_{1},...,X_{n}\right] \\
&  =\mathbb{E}\left[  \mathbb{E}\left\langle \phi-\widehat{\phi}%
,X_{n+1}\right\rangle _{W}^{2}|\widehat{\phi}\right] \\
&  =\mathbb{E}\left[  \left\Vert \Gamma^{1/2}\left(  \phi-\widehat{\phi
}\right)  \right\Vert _{W}^{2}\right]
\end{align*}
Similarly,%
\[
\mathbb{E}\left\langle \psi-\widehat{\psi},X_{n+1}^{\prime}\right\rangle
_{L}^{2}=\mathbb{E}\left[  \left\Vert \Gamma^{\prime\prime1/2}\left(
\psi-\widehat{\psi}\right)  \right\Vert _{L}^{2}\right]
\]
Both preceding equations feature similar expressions. We focus on the term
involving $\phi$ ; we will prove that :%
\[
\mathbb{E}\left[  \left\Vert \Gamma^{1/2}\left(  \phi-\widehat{\phi}\right)
\right\Vert _{W}^{2}\right]  =O\left(  \frac{\beta^{2}}{\alpha^{2}}\right)
+O\left(  \frac{1}{\alpha^{2}\beta^{2}n}\right)  .
\]
Within the proof the reader will easily be convinced that the method would
lead to an analogous result for the term with $\psi.$ From now in order to
alleviate notations we drop the index $\phi$ in $S_{n,\phi}$ and $u_{n,\phi}$.
The sequences $\left(  \alpha_{n}\right)  _{n\in\mathbb{N}}$ and $\left(
\beta_{n}\right)  _{n\in\mathbb{N}}$ will be denoted $\alpha$ and $\beta$
respectively and for short.

We start from%
\begin{align*}
\widehat{\phi}_{n}  &  =\left(  \Gamma_{n}-\Gamma_{n}^{\prime}\Gamma
_{n}^{\prime\prime\dag}\Gamma_{n}^{\prime\ast}+\beta I\right)  ^{-1}\left(
\delta_{n}-\Gamma_{n}^{\prime}\Gamma_{n}^{\prime\prime\dag}\delta_{n}^{\prime
}\right) \\
&  =\left(  S_{n}+\beta I\right)  ^{-1}u_{n}\\
\phi &  =\left(  \Gamma-\Gamma^{\prime}\Gamma^{\prime\prime\dag}\Gamma
^{\prime\ast}\right)  ^{-1}\left(  \delta-\Gamma^{\prime}\Gamma^{\prime
\prime\dag}\delta^{\prime}\right) \\
&  =S^{-1}u
\end{align*}

where we recall that :%
\begin{align*}
u  &  =\delta-\Gamma^{\prime}\Gamma^{\prime\prime\dag}\delta^{\prime},\\
u_{n}  &  =\delta_{n}-\Gamma_{n}^{\prime}\Gamma_{n}^{\prime\prime\dag}%
\delta_{n}^{\prime},\\
S  &  =\Gamma-\Gamma^{\prime}\Gamma^{\prime\prime\dag}\Gamma^{\prime\ast},\\
S_{n}  &  =\Gamma_{n}-\Gamma_{n}^{\prime}\Gamma_{n}^{\prime\prime\dag}%
\Gamma_{n}^{\prime\ast}.
\end{align*}

The proof relies on the following decomposition :%
\begin{align}
\widehat{\phi}_{n}-\phi &  =\left(  S_{n}+\beta I\right)  ^{-1}\left(
u_{n}-u\right)  +\left(  \left(  S_{n}+\beta I\right)  ^{-1}-S^{-1}\right)
u\nonumber\\
&  =\left(  S_{n}+\beta I\right)  ^{-1}\left(  u_{n}-u\right)  +\left(
S_{n}+\beta I\right)  ^{-1}\left(  S-S_{n}-\beta I\right)  S^{-1}u\nonumber\\
&  =A_{n}+B_{n}+C_{n} \label{decomp}%
\end{align}
where%
\begin{align}
A_{n}  &  =\left(  S_{n}+\beta I\right)  ^{-1}\left(  u_{n}-u\right)
\label{T1}\\
B_{n}  &  =\left(  S_{n}+\beta I\right)  ^{-1}\left(  S-S_{n}\right)
\phi\label{T2}\\
C_{n}  &  =\beta\left(  S_{n}+\beta I\right)  ^{-1}\phi\label{T3}%
\end{align}
Along the forthcoming Lemmas we determine rates of convergence for these three
terms. We will prove that the rate of decrease to zero in $L^{2}$ norm is
$\left(  \alpha\beta\sqrt{n}\right)  ^{-1}$ for $A_{n}$ and $B_{n}.$ The rest
of the proof of the main Theorem is postponed to the end of the next and last subsection.

\subsection{Proof of the main Theorem}

The first Lemma gives a reta of convergence for $S_{n}-S$.

\begin{lemma}
\label{cran}The following holds :%
\[
S_{n}-S=\Gamma_{n}-\Gamma_{n}^{\prime}\Gamma_{n}^{\prime\prime\dag}\Gamma
_{n}^{\prime\ast}-\Gamma+\Gamma^{\prime}\Gamma^{\prime\prime\dag}%
\Gamma^{\prime\ast}=O_{L^{2}}\left(  \dfrac{1}{\alpha\sqrt{n}}\right)
\]

\end{lemma}

\begin{proof}
First of all by (\ref{L2}) :%
\[
\left\Vert \Gamma_{n}-\Gamma\right\Vert _{L^{2}\left(  \mathcal{HS}\right)
}=O\left(  \dfrac{1}{\sqrt{n}}\right)
\]
We focus on%
\begin{align*}
&  \Gamma_{n}^{\prime}\Gamma_{n}^{\prime\prime\dag}\Gamma_{n}^{\prime\ast
}-\Gamma^{\prime}\Gamma^{\prime\prime\dag}\Gamma^{\prime\ast}\\
&  =\Gamma_{n}^{\prime}\Gamma_{n}^{\prime\prime\dag}\Gamma_{n}^{\prime\ast
}-\Gamma^{\prime}\Gamma_{n}^{\prime\prime\dag}\Gamma_{n}^{\prime\ast}%
+\Gamma^{\prime}\Gamma_{n}^{\prime\prime\dag}\Gamma_{n}^{\prime\ast}%
-\Gamma^{\prime}\Gamma_{n}^{\prime\prime\dag}\Gamma^{\prime\ast}\\
&  +\Gamma^{\prime}\Gamma_{n}^{\prime\prime\dag}\Gamma^{\prime\ast}%
-\Gamma^{\prime}\Gamma^{\prime\prime\dag}\Gamma^{\prime\ast}.
\end{align*}
Then dealing with each of these three terms separatedly we get%
\begin{align*}
\left\Vert \Gamma_{n}^{\prime}\Gamma_{n}^{\prime\prime\dag}\Gamma_{n}%
^{\prime\ast}-\Gamma^{\prime}\Gamma_{n}^{\prime\prime\dag}\Gamma_{n}%
^{\prime\ast}\right\Vert _{\infty}  &  \leq\left\Vert \Gamma_{n}^{\prime
}-\Gamma^{\prime}\right\Vert _{\infty}\left\Vert \Gamma_{n}^{\prime\prime\dag
}\Gamma_{n}^{\prime\ast}\right\Vert _{\infty}\\
&  \leq\left\Vert \Gamma_{n}^{\prime}-\Gamma^{\prime}\right\Vert _{\infty
}\left\Vert \left(  \Gamma_{n}^{\prime\prime\dag}\right)  \right\Vert
_{\infty}\left\Vert \Gamma_{n}^{\prime\ast}\right\Vert _{\infty}\\
&  \leq C\dfrac{\left\Vert \Gamma_{n}^{\prime}-\Gamma^{\prime}\right\Vert
_{\infty}}{\alpha}\quad a.s.
\end{align*}
The last bound was derived from (\ref{norm.delta}) and (\ref{norm.res}).%
\begin{align*}
&  \left\Vert \Gamma^{\prime}\Gamma_{n}^{\prime\prime\dag}\left(  \Gamma
_{n}^{\prime\ast}-\Gamma^{\prime\ast}\right)  \right\Vert _{\infty}\\
&  \leq C\dfrac{\left\Vert \Gamma_{n}^{\prime\ast}-\Gamma^{\prime\ast
}\right\Vert _{\infty}}{\alpha}\quad a.s.
\end{align*}
At last,%
\begin{align*}
\Gamma^{\prime}\left(  \Gamma_{n}^{\prime\prime\dag}-\Gamma^{\prime\prime\dag
}\right)  \Gamma^{\prime\ast}  &  =\Gamma^{\prime}\Gamma_{n}^{\prime\prime
\dag}\left(  \Gamma^{\prime\prime}-\Gamma_{n}^{\prime\prime}\right)
\Gamma^{\prime\prime\dag}\Gamma^{\prime\ast}\\
&  =\Gamma^{\prime}\Gamma_{n}^{\prime\prime\dag}\left(  \Gamma^{\prime\prime
}-\Gamma_{n}^{\prime\prime}\right)  \Gamma^{\prime\prime\dag}\Gamma
^{\prime\ast}\\
&  =\left(  \Gamma^{\prime}-\Gamma_{n}^{\prime}\right)  \Gamma_{n}%
^{\prime\prime\dag}\left(  \Gamma^{\prime\prime}-\Gamma_{n}^{\prime\prime
}\right)  \Gamma^{\prime\prime\dag}\Gamma^{\prime\ast}\\
&  +\Gamma_{n}^{\prime}\Gamma_{n}^{\prime\prime\dag}\left(  \Gamma
^{\prime\prime}-\Gamma_{n}^{\prime\prime}\right)  \Gamma^{\prime\prime\dag
}\Gamma^{\prime\ast}%
\end{align*}
Then,%
\begin{align*}
&  \left\Vert \Gamma^{\prime}\left(  \Gamma_{n}^{\prime\prime\dag}%
-\Gamma^{\prime\prime\dag}\right)  \Gamma^{\prime\ast}\right\Vert _{\infty}\\
&  \leq\left\Vert \left(  \Gamma^{\prime}-\Gamma_{n}^{\prime}\right)
\Gamma_{n}^{\prime\prime\dag}\left(  \Gamma^{\prime\prime}-\Gamma_{n}%
^{\prime\prime}\right)  \Gamma^{\prime\prime\dag}\Gamma^{\prime\ast
}\right\Vert _{\infty}\\
&  +\left\Vert \Gamma_{n}^{\prime}\left(  \Gamma_{n}^{\prime\prime\dag
}\right)  ^{1/2}\right\Vert _{\infty}\left\Vert \left(  \Gamma_{n}%
^{\prime\prime\dag}\right)  ^{1/2}\left(  \Gamma^{\prime\prime}-\Gamma
_{n}^{\prime\prime}\right)  \left(  \Gamma^{\prime\prime\dag}\right)
^{1/2}\right\Vert _{\infty}\left\Vert \left(  \Gamma^{\prime\prime\dag
}\right)  ^{1/2}\Gamma^{\prime\ast}\right\Vert _{\infty}%
\end{align*}
By Proposition \ref{lodge} the second term may be bounded by%
\[
C\left\Vert \left(  \Gamma_{n}^{\prime\prime\dag}\right)  ^{1/2}\left(
\Gamma^{\prime\prime}-\Gamma_{n}^{\prime\prime}\right)  \left(  \Gamma
^{\prime\prime\dag}\right)  ^{1/2}\right\Vert _{\infty}=O_{L^{2}}\left(
\dfrac{1}{\alpha\sqrt{n}}\right)
\]
since
\[
\left\Vert \Gamma_{n}^{\prime\prime\dag1/2}\right\Vert _{\infty}=\left\Vert
\Gamma^{\prime\prime\dag1/2}\right\Vert _{\infty}=\alpha^{-1/2}.
\]
Cauchy-Schwartz inequality yields for the first :%
\begin{align*}
&  \mathbb{E}\left\Vert \left(  \Gamma^{\prime}-\Gamma_{n}^{\prime}\right)
\Gamma_{n}^{\prime\prime\dag}\left(  \Gamma^{\prime\prime}-\Gamma_{n}%
^{\prime\prime}\right)  \Gamma^{\prime\prime\dag}\Gamma^{\prime\ast
}\right\Vert _{\infty}^{2}\\
&  \leq M\left(  \mathbb{E}\left\Vert \left(  \Gamma^{\prime}-\Gamma
_{n}^{\prime}\right)  \Gamma_{n}^{\prime\prime\dag}\right\Vert _{\infty}%
^{4}\mathbb{E}\left\Vert \left(  \Gamma^{\prime\prime}-\Gamma_{n}%
^{\prime\prime}\right)  \Gamma^{\prime\prime\dag}\right\Vert _{\infty}%
^{4}\right)  ^{1/2}\\
&  \leq M\dfrac{1}{n^{2}\alpha^{4}},
\end{align*}
hence%
\[
\left\Vert \left(  \Gamma^{\prime}-\Gamma_{n}^{\prime}\right)  \Gamma
_{n}^{\prime\prime\dag}\left(  \Gamma^{\prime\prime}-\Gamma_{n}^{\prime\prime
}\right)  \Gamma^{\prime\prime\dag}\Gamma^{\prime\ast}\right\Vert =O_{L^{2}%
}\left(  \dfrac{1}{\alpha^{2}n}\right)  .
\]
The proof of Lemma \ref{cran} is finished.
\end{proof}

\begin{lemma}
\label{u}We have :%
\[
u_{n}-u=O_{L^{2}}\left(  \dfrac{1}{\alpha\sqrt{n}}\right)  .
\]

\end{lemma}

\begin{proof}
We start with :%
\[
u_{n}-u=\delta_{n}-\delta+\Gamma^{\prime}\Gamma^{\prime\prime\dag}%
\delta^{\prime}-\Gamma_{n}^{\prime}\Gamma_{n}^{\prime\prime\dag}\delta
_{n}^{\prime}.
\]
Clearly $\delta_{n}-\delta=O_{L^{2}}\left(  \dfrac{1}{\sqrt{n}}\right)  $ and
we study the second term%
\begin{align*}
\Gamma_{n}^{\prime}\Gamma_{n}^{\prime\prime\dag}\delta_{n}^{\prime}%
-\Gamma^{\prime}\Gamma^{\prime\prime\dag}\delta^{\prime}  &  =\left(
\Gamma_{n}^{\prime}-\Gamma^{\prime}\right)  \Gamma_{n}^{\prime\prime\dag
}\delta_{n}^{\prime}+\Gamma^{\prime}\left(  \Gamma_{n}^{\prime\prime\dag
}-\Gamma^{\prime\prime\dag}\right)  \delta_{n}^{\prime}\\
&  +\Gamma^{\prime}\Gamma^{\prime\prime\dag}\left(  \delta_{n}^{\prime}%
-\delta^{\prime}\right)  .
\end{align*}
Since $\delta_{n}^{\prime}$ is almost surely bounded (see (\ref{asd})),
$\Gamma_{n}^{\prime}-\Gamma=O_{L^{2}}\left(  \dfrac{1}{\sqrt{n}}\right)  ,$
$\delta_{n}^{\prime}-\delta^{\prime}=O_{L^{2}}\left(  \dfrac{1}{\sqrt{n}%
}\right)  $ and $\left\Vert \Gamma^{\prime\prime\dag}\right\Vert _{\infty
}=\left\Vert \Gamma_{n}^{\prime\prime\dag}\right\Vert _{\infty}=\alpha^{-1}$
we get :%
\begin{align*}
\left\Vert \left(  \Gamma_{n}^{\prime}-\Gamma^{\prime}\right)  \Gamma
_{n}^{\prime\prime\dag}\delta_{n}^{\prime}\right\Vert _{W}  &  =O_{L^{2}%
}\left(  \dfrac{1}{\alpha\sqrt{n}}\right)  ,\\
\left\Vert \Gamma^{\prime}\Gamma^{\prime\prime\dag}\left(  \delta_{n}^{\prime
}-\delta^{\prime}\right)  \right\Vert _{W}  &  =O_{L^{2}}\left(  \dfrac
{1}{\alpha\sqrt{n}}\right)  .
\end{align*}
The remaining term is%
\begin{align*}
\Gamma^{\prime}\left(  \Gamma_{n}^{\prime\prime\dag}-\Gamma^{\prime\prime\dag
}\right)  \delta_{n}^{\prime}  &  =\Gamma^{\prime}\Gamma^{\prime\prime\dag
}\left(  \Gamma^{\prime\prime}-\Gamma_{n}^{\prime\prime}\right)  \Gamma
_{n}^{\prime\prime\dag}\delta_{n}^{\prime}\\
&  =\Gamma^{\prime}\Gamma^{\prime\prime\dag}\left(  \Gamma^{\prime\prime
}-\Gamma_{n}^{\prime\prime}\right)  \Gamma_{n}^{\prime\prime\dag}\left(
\Gamma_{n}^{\prime\ast}\phi+\Gamma_{n}^{\prime\prime}\psi+u_{n}^{\prime
}\right)  ,\\
&  =\Gamma^{\prime}\left(  \Gamma^{\prime\prime\dag}\right)  ^{1/2}\left(
m_{1}+m_{2}+m_{3}\right)
\end{align*}
where
\begin{align*}
m_{1}  &  =\left(  \Gamma^{\prime\prime\dag}\right)  ^{1/2}\left(
\Gamma^{\prime\prime}-\Gamma_{n}^{\prime\prime}\right)  \left(  \Gamma
_{n}^{\prime\prime\dag}\right)  ^{1/2}\left(  \Gamma_{n}^{\prime\prime\dag
}\right)  ^{1/2}\Gamma_{n}^{\prime\ast}\phi,\\
m_{2}  &  =\left(  \Gamma^{\prime\prime\dag}\right)  ^{1/2}\left(
\Gamma^{\prime\prime}-\Gamma_{n}^{\prime\prime}\right)  \Gamma_{n}%
^{\prime\prime\dag}\Gamma_{n}^{\prime\prime}\psi,\\
m_{3}  &  =\left(  \Gamma^{\prime\prime\dag}\right)  ^{1/2}\left(
\Gamma^{\prime\prime}-\Gamma_{n}^{\prime\prime}\right)  \Gamma_{n}%
^{\prime\prime\dag}e_{n}^{\prime}.
\end{align*}
First we drop $\Gamma^{\prime}\left(  \Gamma^{\prime\prime\dag}\right)
^{1/2}$ since the norm of this operator may be bounded by a constant
independent from $\alpha$ (see Proposition \ref{lodge}). We turn to :%
\begin{align*}
\left\Vert m_{1}\right\Vert  &  \leq M\left\Vert \left(  \Gamma^{\prime\prime
}-\Gamma_{n}^{\prime\prime}\right)  \right\Vert _{\infty}\left\Vert \left(
\Gamma_{n}^{\prime\prime\dag}\right)  ^{1/2}\right\Vert _{\infty}\left\Vert
\left(  \Gamma^{\prime\prime\dag}\right)  ^{1/2}\right\Vert _{\infty},\\
\left\Vert m_{2}\right\Vert  &  \leq\left\Vert \psi\right\Vert _{L}\left\Vert
\left(  \Gamma^{\prime\prime\dag}\right)  ^{1/2}\right\Vert _{\infty
}\left\Vert \left(  \Gamma^{\prime\prime}-\Gamma_{n}^{\prime\prime}\right)
\right\Vert _{\infty}%
\end{align*}
since $\left\Vert \Gamma_{n}^{\prime\prime\dag}\Gamma_{n}^{\prime\prime
}\right\Vert _{\infty}\leq1$ almost surely. The consequence of the display
above is $\left\Vert m_{1}\right\Vert _{L^{2}}=O\left(  \dfrac{1}{\alpha
\sqrt{n}}\right)  $ and $\left\Vert m_{2}\right\Vert _{L^{2}}=O\left(
\dfrac{1}{\sqrt{\alpha n}}\right)  .$

We can deal with $m_{3}$ as was done within the proof of the preceding Lemma
\ref{cran}. Clearly we may cope with $m_{3}$ as if the random $\Gamma
_{n}^{\prime\prime\dag}$ was replaced by the non random $\Gamma^{\prime
\prime\dag}$. We should study%
\[
\left[  \left(  \Gamma^{\prime\prime\dag}\right)  ^{1/2}\left(  \Gamma
^{\prime\prime}-\Gamma_{n}^{\prime\prime}\right)  \left(  \Gamma^{\prime
\prime\dag}\right)  ^{1/2}\right]  \left[  \left(  \Gamma^{\prime\prime\dag
}\right)  ^{1/2}e_{n}^{\prime}\right]  .
\]
It is enough to get a rate of decrease for each of the these terms. Once again
we have :%
\begin{align*}
\left\Vert \left(  \Gamma^{\prime\prime\dag}\right)  ^{1/2}\left(
\Gamma^{\prime\prime}-\Gamma_{n}^{\prime\prime}\right)  \left(  \Gamma
^{\prime\prime\dag}\right)  ^{1/2}\right\Vert _{L^{2}}  &  =O\left(  \frac
{1}{\alpha\sqrt{n}}\right) \\
\left\Vert \left(  \Gamma^{\prime\prime\dag}\right)  ^{1/2}e_{n}^{\prime
}\right\Vert _{L^{2}}  &  =O\left(  \frac{1}{\sqrt{\alpha n}}\right)
\end{align*}
which completes the proof of Lemma \ref{u}.
\end{proof}

Now we are ready to go back to (\ref{T1}) and (\ref{T2}) as announced sooner.

\begin{lemma}
\label{A_n}We have :%
\begin{align*}
A_{n}  &  =O_{L^{2}}\left(  \frac{1}{\alpha\beta\sqrt{n}}\right)  ,\\
B_{n}  &  =O_{L^{2}}\left(  \frac{1}{\alpha\beta\sqrt{n}}\right)  .
\end{align*}

\end{lemma}

\begin{proof}
Since
\[
\left\Vert A_{n}\right\Vert \leq\left\Vert \left(  S_{n}+\beta I\right)
^{-1}\right\Vert _{\infty}\left\Vert u_{n}-u\right\Vert _{W}%
\]
by Lemma \ref{u} and Proposition \ref{KT} we get the first desired result
\end{proof}

Once again the proof of the second relies on Proposition \ref{KT} and Lemma
\ref{cran}. Indeed%
\begin{align*}
\left\Vert B_{n}\right\Vert _{W}  &  \leq\left\Vert \left(  S_{n}+\beta
I\right)  ^{-1}\right\Vert _{\infty}\left\Vert S-S_{n}\right\Vert _{\infty
}\left\Vert \phi\right\Vert _{W}\\
&  \leq\dfrac{\left\Vert \phi\right\Vert _{W}}{\beta}\left\Vert S-S_{n}%
\right\Vert _{\infty}%
\end{align*}
hence the result.

We should deal with the last term. In a first step we prove that $S_{n}$ may
be replaced by $S.$

\begin{lemma}
\label{C_n}When $\alpha\beta\sqrt{n}\rightarrow+\infty$,%
\[
C_{n}=\beta\left(  S+\beta I\right)  ^{-1}\phi\left(  1+o\left(  1\right)
\right)  .
\]

\end{lemma}

\begin{remark}
The preceding equality should be understood with respect to the $L^{2}$ norm.
\end{remark}

\begin{proof}
Successively,%
\begin{align*}
C_{n}  &  =\beta\left(  S_{n}+\beta I\right)  ^{-1}\phi\\
&  =\beta\left(  \left(  S_{n}+\beta I\right)  ^{-1}-\left(  S+\beta I\right)
^{-1}\right)  \phi+\beta\left(  S+\beta I\right)  ^{-1}\phi\\
&  =\left[  \left(  \left(  S_{n}+\beta I\right)  ^{-1}\left(  S-S_{n}\right)
\right)  +I\right]  \beta\left(  S+\beta I\right)  ^{-1}\phi
\end{align*}
and
\[
\left\Vert C_{n}\right\Vert \leq\left\Vert \beta\left(  S+\beta I\right)
^{-1}\phi\right\Vert _{W}\left(  1+\left\Vert \left(  S_{n}+\beta I\right)
^{-1}\left(  S-S_{n}\right)  \right\Vert _{\infty}\right)  .
\]
Now it suffices to apply Lemma \ref{A_n} to get the desired result.
\end{proof}

The next Lemma may be hard to understand at first glance. Within the
forthcoming proof of Theorem \ref{conv} the bias term $C_{n}$ will slightly
change. We refer to displays (\ref{E1}) and (\ref{E2}) below for a deeper understanding.

\begin{lemma}
\label{placebo}The following holds :%
\[
\left\Vert \Gamma^{1/2}\left(  S+\beta I\right)  ^{-1}\Gamma^{1/2}\right\Vert
_{\infty}=O\left(  \frac{1}{\alpha}\right)  .
\]

\end{lemma}

\begin{proof}
Once again it takes two steps to get the result. First note that $\Gamma
^{1/2}\left(  S\right)  ^{-1}\Gamma^{1/2}$ is a bounded linear operator.
Indeed%
\begin{equation}
S=\Gamma-\Gamma^{\prime}\Gamma^{\prime\prime\dag}\Gamma^{\prime\ast}%
=\Gamma^{1/2}\Lambda_{\alpha}\Gamma^{1/2} \label{queer}%
\end{equation}
where $R=D\Gamma^{1/2},$%
\[
\Lambda_{\alpha}=I-R^{\ast}\left(  RR^{\ast}+\alpha I\right)  ^{-1}R.
\]
The Schmidt decomposition of $R$ is (see (\ref{mano}) above for the empirical
version) :%
\[
R=\sum_{k\in\mathbb{N}}\mu_{k}\left(  u_{k}\otimes v_{k}\right)  .
\]
where $\left(  u_{k}\right)  _{k\in\mathbb{N}}$ (resp. $\left(  v_{k}\right)
_{k\in\mathbb{N}}$) is a complete orthonormal system in $W$ (resp. $L$). Hence
:%
\begin{align*}
\Lambda_{\alpha}  &  =\sum_{k\in\mathbb{N}}\left(  1-\dfrac{\mu_{k}^{2}}%
{\mu_{k}^{2}+\alpha}\right)  \left(  u_{k}\otimes u_{k}\right) \\
&  =\sum_{k\in\mathbb{N}}\dfrac{\alpha}{\mu_{k}^{2}+\alpha}\left(
u_{k}\otimes u_{k}\right)  .
\end{align*}
The operator $\Lambda_{\alpha}$ has a bounded inverse%
\[
\Lambda_{\alpha}^{-1}=\dfrac{1}{\alpha}\sum_{k\in\mathbb{N}}\left(  \mu
_{k}^{2}+\alpha\right)  \left(  u_{k}\otimes u_{k}\right)
\]
and $\left\Vert \Lambda_{\alpha}^{-1}\right\Vert _{\infty}=1+\left(  \sup
\mu_{k}^{2}\right)  /\alpha\leq M/\alpha$ for $M$ large enough (or $\alpha$
small enough).\newline Hence%
\begin{equation}
\Gamma^{1/2}\left(  S\right)  ^{-1}\Gamma^{1/2}=\Gamma^{1/2}\Gamma
^{-1/2}\Lambda_{\alpha}^{-1}\Gamma^{-1/2}\Gamma^{1/2}=\Lambda_{\alpha}^{-1}.
\label{BH}%
\end{equation}
Now (second step) we prove that :%
\[
\Gamma^{1/2}\left(  S+\beta I\right)  ^{-1}\Gamma^{1/2}\ll\Gamma^{1/2}%
S^{-1}\Gamma^{1/2}.
\]
Let us pick a given $x$ in $W$, then
\[
\left\langle \Gamma^{1/2}\left(  S+\beta I\right)  ^{-1}\Gamma^{1/2}%
x,x\right\rangle _{W}=\left\langle \left(  S+\beta I\right)  ^{-1}\Gamma
^{1/2}x,\Gamma^{1/2}x\right\rangle _{W}%
\]
It suffices to get for all $y$ in in the domain of operator $\Gamma^{-1/2}$ :%
\begin{equation}
\left\langle \left(  S+\beta I\right)  ^{-1}y,y\right\rangle _{W}%
\leq\left\langle S^{-1}y,y\right\rangle _{W} \label{CT}%
\end{equation}
Standard results on the spectrum of $\left(  S+\beta I\right)  ^{-1}S$ prove
that $\left(  S+\beta I\right)  ^{-1}S\geq0$ and that $\left\Vert \left(
S+\beta I\right)  ^{-1}S\right\Vert \leq1$ which is enough to claim
(\ref{CT}).\newline We are now in position to finixh the proof of the Lemma.
It is plain from (\ref{CT}) that%
\[
\left\Vert \Gamma^{1/2}\left(  S+\beta I\right)  ^{-1}\Gamma^{1/2}\right\Vert
_{\infty}\leq\left\Vert \Gamma^{1/2}\left(  S\right)  ^{-1}\Gamma
^{1/2}\right\Vert _{\infty}=\left\Vert \Lambda_{\alpha}^{-1}\right\Vert
_{\infty}\leq\frac{C}{\alpha}%
\]
which is the claimed result.\bigskip
\end{proof}

\textbf{Proof of Theorem }\ref{conv}\textbf{:}

Now starting from (\ref{decomp}) we get%
\begin{align}
\left\Vert \Gamma^{1/2}\left(  \phi-\widehat{\phi}\right)  \right\Vert
_{W}^{2}  &  \leq M\left\Vert \Gamma^{1/2}\left(  A_{n}+B_{n}+C_{n}\right)
\right\Vert _{W}^{2}\nonumber\\
&  \leq M\left(  \left\Vert A_{n}\right\Vert _{W}^{2}+\left\Vert
B_{n}\right\Vert _{W}^{2}+\left\Vert \Gamma^{1/2}C_{n}\right\Vert _{W}%
^{2}\right)  . \label{E1}%
\end{align}

Lemmas \ref{A_n} gives the rates of convergence for $\left\Vert A_{n}%
\right\Vert _{W}^{2}$ and $\left\Vert B_{n}\right\Vert _{W}^{2}$ respectively.
But Lemma \ref{C_n} is unfortunately not enough to get a rate in the last
term. However this previous Lemma enables to focus on :%
\begin{equation}
\beta\Gamma^{1/2}\left(  S+\beta I\right)  ^{-1}\phi=\beta\Gamma^{1/2}\left(
S+\beta I\right)  ^{-1}\Gamma^{1/2}\Gamma^{-1/2}\phi\label{E2}%
\end{equation}
and
\begin{equation}
\left\Vert \Gamma^{1/2}C_{n}\right\Vert _{W}^{2}\leq M\beta^{2}\left\Vert
\Gamma^{1/2}\left(  S+\beta I\right)  ^{-1}\Gamma^{1/2}\right\Vert _{\infty
}^{2}\left\Vert \Gamma^{-1/2}\phi\right\Vert _{W}. \label{marin}%
\end{equation}
By assumption $\mathbf{A4}$, $\left\Vert \Gamma^{-1/2}\phi\right\Vert _{W}$ is
finite. We deal with the central term, namely :%
\begin{align*}
\Gamma^{1/2}\left(  S+\beta I\right)  ^{-1}\Gamma^{1/2}  &  =\Gamma
^{1/2}\left(  \Gamma^{1/2}\Lambda_{\alpha}\Gamma^{1/2}+\beta I\right)
^{-1}\Gamma^{1/2}\\
&  \ll\Gamma^{1/2}\left(  \Gamma^{1/2}\Lambda_{\alpha}\Gamma^{1/2}\right)
^{-1}\Gamma^{1/2}=\Lambda_{\alpha}^{-1}.
\end{align*}
(see (\ref{queer})) and%
\[
\left\Vert \Lambda_{\alpha}^{-1}\right\Vert _{\infty}^{2}=O\left(  \alpha
^{-2}\right)  .
\]
Collecting this last display with (\ref{marin}) we get%
\[
\left\Vert \Gamma^{1/2}C_{n}\right\Vert _{W}^{2}=O\left(  \frac{\beta^{2}%
}{\alpha^{2}}\right)  .
\]
This finishes the proof of Theorem \ref{conv}.


\begin{thebibliography}{99}                                                                                               %


\bibitem {AF}Adams R.A. and Fournier J.J.F., 2003. Sobolev spaces, Academic
Press, 2nd ed.

\bibitem {Amato.et.al.CSDA}Amato U., Antoniadis A., Feiss I., 2006. Dimension
reduction in functional regression with applications, to appear in
\textit{Computational Statistics and Data Analysis.}

\bibitem {Bos3}Bosq, D., 2000. \textit{Linear processes in function spaces}.
Lectures notes in statistics. Springer Verlag.

\bibitem {Brown.Fearn.Vanucci}Brown, P.J., Fearn, T. and Vanucci M., 2001.
Bayesian wavelet regression on curves with application to a spectroscopic
calibration problem, \textit{Journal of the American Statistical
Association.}, 96 (454), 398--408.

\bibitem {CH}Cai T., Hall, P., 2006. Prediction in functional linear
regression, \textit{Annals of Statistics}, \textbf{34}, n$%
{{}^\circ}%
5$.

\bibitem {Cardot.SEFL}Cardot H., Ferraty F., Sarda P., 2006. Spline estimators
for the functional linear model: Consistency, Application and Splus
implementation, To appear in \textit{Statistica Sinica}.

\bibitem {CMS2006}Cardot H., Mas A., Sarda P, 2006. Weak convergence in the
functional linear model. To appear in \textit{Probability Theory and Related
Fields}.

\bibitem {DS}Dunford, N. and Schwartz, J.T., 1988. \textit{Linear Operators,
Vol. I \& II}. Wiley Classics Library.

\bibitem {FV}Ferraty, F., Vieu P., 2003. The functional nonparametric model
and application to spectrometric data. \textit{Computational Statistics}
\textbf{17} n$%
{{}^\circ}%
4,$ 545-564.

\bibitem {FV2}Ferraty, F., Vieu P., 2006. Nonparametic Functional Data
Analysis, Springer.

\bibitem {Frank.Friedman}Frank I.E. and Friedman J.H., 1993. A statistical
view of some chemometrics regression tools, \textit{Technometrics},
\textbf{35}, no 2, 109-135.

\bibitem {GGK}Gohberg, I., Goldberg, S. and Kaashoek,M.A., 1991.
\textit{Classes of linear operators Vol I \& II. Operator Theory : advances
and applications}, Birkha\"{u}ser Verlag.

\bibitem {Hastie.Mallows}Hastie T. and Mallows C., 1993. Discussion on the
paper of Frank I.E. and Friedman J.H., 1993, A statistical view of some
chemometrics regression tools, \textit{Technometrics}, \textbf{35}, no 2, 140-143.

\bibitem {LT}Ledoux, M. and Talagrand, M, 1991 : \textit{Probability in Banach
spaces. Isoperimetry and processes}. Ergebnisse der Mathematik und ihrer
Grenzgebiete, \textbf{23}. Springer-Verlag, Berlin.

\bibitem {Marx.Eilers}Marx, B.D. and Eilers, P.H.C. (2002) Multivariate
calibration stability: a comparison of methods, \textit{Journal of
Chemometrics}, \textbf{16}, 129-140.

\bibitem {MP}Mas A, Pumo, B, 2006. The ARHD model. To appear in JSPI.
Technical report available at
\textit{http://fr.arxiv.org/PS\_cache/math/pdf/0502/0502285.pdf}

\bibitem {marion.pumo}Marion J.M., Pumo B., 2004. Comparaison des mod\`{e}les
ARH(1) et ARHD(1) sur des donn\'{e}es physiologiques, \textit{Annales de
l'ISUP}, \textbf{48}, 3, pp. 29-38.

\bibitem {Osborne.et.al.}Osborne, B.J., Fearn, T., Miller, A.R. and Douglas,
S. (1984) Application of Near-Infrared Reflectance Spectroscopy to
Compositional Analysis of Biscuits Dougts, \textit{J. of the Sc. of Food and
Agricult.}, \textbf{35}, 99--105

\bibitem {Ram1}Ramsay J.O., 2000. Differential equation models for statistical
functions, \textit{Canadian Journal of Statistics}, \textbf{28}, n$%
{{}^\circ}%
2,$ 225-240.

\bibitem {RamDal}Ramsay J.O., Dalzell C.J., 1991, Some tools for functional
data analysis (with dsicussion), \textit{Journal of the Royal Statistical
Society}, B, \textbf{53}, 539-572.

\bibitem {ramsey}Ramsay J.O., Silverman B.W., 1997. \textit{Functional Data
Analysis}, Springer.

\bibitem {ramsay2}Ramsay J.O., Silverman B.W., 2002. \textit{Applied
Functional Data Analysis: Methods and Case Studies}, Springer.

\bibitem {Silv}Silverman B.W., 1996. Smoothed functional principal component
analysis by choice of norm, \textit{Annals of Statistics}, \textbf{24}, No 1, 1-24.

\bibitem {TA}Tikhonov A.N., Arsenin V.Y., 1977. \textit{Solutions of ill-posed
problems}. V.H. Winstons and sons, Washington..

\bibitem {Z}Ziemer W.P., 1989. Weakly differentiable functions.
\textit{Sobolev spaces and functions of bounded variations}. Graduate Text in
Mathematics 120. Springer-Verlag, New-York.
\end{thebibliography}
\end{document}